# Pascal's triangle, Hoggatt matrices, and analogous constructions

Johann Cigler

**Abstract**

We give an overview about some elementary properties of Hoggatt matrices, which are generalizations of Pascal's triangle, and study $q-$analogs and Fibonacci analogs and derive a common generalization.

## 1. Introduction

In [4] Daniel C. Fielder and Cecil O. Alford defined generalizations of Pascal's triangle which they called *Hoggatt triangles*. We give an overview about some elementary properties of these triangles and their $q-$analogs and give a common generalization with Fibonacci polynomials. I want to thank Christian Krattenthaler for help with some determinants and hypergeometric identities.

Let us first introduce some notations which emphasize the analogy with Pascal's triangle. Let $d$ be a positive integer. We write

$$\langle n\rangle_d=\binom{n+d-1}{d}=\frac{(d+1)(d+2)\cdots(d+n-1)}{(n-1)!} \tag{1}$$

and

$$\langle n\rangle_d!=\prod_{j=1}^{n}\langle j\rangle_d \tag{2}$$

and define

$$\left\langle {n \atop k} \right\rangle_d=\frac{\langle n\rangle_d}{\langle k\rangle_d}\left\langle {n-1 \atop k-1} \right\rangle_d=\prod_{j=0}^{k-1}\frac{\langle n-j\rangle_d}{\langle k-j\rangle_d}=\frac{\langle n\rangle_d!}{\langle k\rangle_d!\langle n-k\rangle_d!} \tag{3}$$

for $0\le k\le n$ and $\left\langle {n \atop k} \right\rangle_d=0$ for $k>n.$

Following [4] we call the matrix

$$H_d=\left(\left\langle {n \atop k} \right\rangle_d\right)_{n,k\ge 0}$$

the *Hoggatt matrix* or *Hoggatt triangle* of order $d$.

--------------------------------

*Email*: johann.cigler@univie.ac.at

All entries of these matrices are nonnegative integers. This is of course true for Pascal's triangle $H_1=\left(\binom{n}{k}\right)_{n,k\geq 0}$ because of the recursion $\binom{n}{k}=\binom{n-1}{k}+\binom{n-1}{k-1}.$

For $d=2$ the numbers $\langle n\rangle_2=\binom{n+1}{2}=T_n$ are the triangle numbers $1,3,6,10,\cdots$ and $\langle n\rangle_2!=\frac{n!(n+1)!}{2^n}$, because $\langle 1\rangle_2!=1$ and by induction

$$\langle n\rangle_2!=\langle n-1\rangle_2!\langle n\rangle_2=\frac{(n-1)!n!}{2^{n-1}}\frac{n(n+1)}{2}=\frac{n!(n+1)!}{2^n}.$$

The entries

$$\left\langle \begin{matrix} n \\ k \end{matrix} \right\rangle_2=\frac{n!(n+1)!}{2^n}\frac{2^k}{k!(k+1)!}\frac{2^{n-k}}{(n-k)!(n-k+1)!}=\frac{1}{k+1}\binom{n}{k}\binom{n+1}{k}$$

are Narayana numbers and $H_2$ is the well-known Narayana triangle (cf. e.g. [9] and OEIS [8], A001263)

$$\begin{matrix} & & & & & 1 & & & & & \\ & & & & 1 & & 1 & & & & \\ & & & 1 & & 3 & & 1 & & & \\ & & 1 & & 6 & & 6 & & 1 & & \\ & 1 & & 10 & & 20 & & 10 & & 1 & \\ 1 & & 15 & & 50 & & 50 & & 15 & & 1 \end{matrix}.$$

For $3\leq d\leq 9$ the corresponding triangles appear in OEIS [8], A056939, A056940, A056941, A142465, A142467, A142468, A174109.

From $\left\langle \begin{matrix} n \\ k \end{matrix} \right\rangle_d=\frac{(d+n-k)\cdots(d+n-1)}{(n-k)\cdots(n-1)}\left\langle \begin{matrix} n-1 \\ k \end{matrix} \right\rangle_d$ we see that the first terms of $H_d$ are

$$\begin{pmatrix} 1 & & & & & \\ 1 & 1 & & & & \\ 1 & d+1 & 1 & & & \\ 1 & \frac{(d+1)(d+2)}{2} & \frac{(d+1)(d+2)}{2} & 1 & & \\ 1 & \frac{(d+1)(d+2)(d+3)}{6} & \frac{(d+1)(d+2)^2(d+3)}{12} & \frac{(d+1)(d+2)(d+3)}{6} & 1 & \\ 1 & \frac{(d+1)(d+2)(d+3)(d+4)}{24} & \frac{(d+1)(d+2)^2(d+3)^2(d+4)}{144} & \frac{(d+1)(d+2)^2(d+3)^2(d+4)}{144} & \frac{(d+1)(d+2)(d+3)(d+4)}{24} & 1 \end{pmatrix}$$

It is also well-known and easy to verify that

$$\left\langle {n \atop k} \right\rangle_2 = \det\begin{pmatrix} \binom{n}{k} & \binom{n+1}{k+1} \\ \binom{n+1}{k} & \binom{n+2}{k+1} \end{pmatrix} = \det\begin{pmatrix} \binom{n}{k} & \binom{n}{k-1} \\ \binom{n}{k+1} & \binom{n}{k} \end{pmatrix} \tag{4}$$

These determinants show that all elements $\left\langle {n \atop k} \right\rangle_2$ are integers.

## 2. Main properties

For each $n$ the entries $\left\langle {n \atop k} \right\rangle_d$ are *palindromic with center of symmetry at* $\frac{n}{2}$ since

$$\left\langle {n \atop k} \right\rangle_d = \left\langle {n \atop n-k} \right\rangle_d . \tag{5}$$

They are also *unimodal* with center of symmetry at $\frac{n}{2}$, which means that

$$\left\langle {n \atop 0} \right\rangle_d \le \left\langle {n \atop 1} \right\rangle_d \le \cdots \le \left\langle {n \atop \left\lfloor \frac{n}{2} \right\rfloor} \right\rangle_d = \left\langle {n \atop \left\lfloor \frac{n+1}{2} \right\rfloor} \right\rangle_d \ge \cdots \ge \left\langle {n \atop n-1} \right\rangle_d \ge \left\langle {n \atop n} \right\rangle_d .$$

Due to symmetry it suffices to show that we have $\left\langle {n \atop k} \right\rangle_d \le \left\langle {n \atop k+1} \right\rangle_d$ or equivalently

$\langle k+1 \rangle_d \le \langle n-k \rangle_d$ for $k < \left\lfloor \frac{n}{2} \right\rfloor$.

This is true because for each $j$ we have $k+1+j \le n-k+j.$

Let us also mention some alternative formulas.

**Proposition 1**

$$\langle n \rangle_d ! = \prod_{j=0}^{d-1} \frac{(n+j)!}{(d-j)^{n+j}}. \tag{6}$$

*and*

$$\left\langle {n \atop k} \right\rangle_d = \prod_{j=0}^{k-1} \frac{\langle n-j \rangle_d}{\langle k-j \rangle_d} = \prod_{j=0}^{d-1} \frac{\binom{n+j}{k}}{\binom{k+j}{k}} = \prod_{j=0}^{d-1} \frac{\binom{n+d-1}{k+j}}{\binom{n+d-1}{j}}. \tag{7}$$

**Proof**

Let $f(n)$ denote the right-hand side of (6). Then $f(1)=\prod_{j=0}^{d-1}\frac{(1+j)!}{(d-j)^{1+j}}=1=\langle 1\rangle_d!$ and by induction

$$f(n)=\prod_{j=0}^{d-1}\frac{(n+j)!}{(d-j)^{n+j}}=\prod_{j=0}^{d-1}\frac{(n+j-1)!}{(d-j)^{n+j-1}}\frac{(n+d-1)!}{(n-1)!d!}=\binom{n+d-1}{d}f(n-1)=\langle n\rangle_d\langle n-1\rangle_d!.$$

Identities (7) follow from

$$\prod_{j=0}^{d-1}\frac{\binom{n+j}{k}}{\binom{k+j}{k}}=\prod_{j=0}^{d-1}\frac{\frac{n+j}{k}\binom{n+j-1}{k-1}}{\frac{k+j}{k}\binom{k+j-1}{k-1}}=\left\langle\begin{matrix}n-1\\k-1\end{matrix}\right\rangle_d\prod_{j=0}^{d-1}\frac{n+j}{k+j}=\left\langle\begin{matrix}n-1\\k-1\end{matrix}\right\rangle_d\frac{\binom{n+d-1}{d}}{\binom{k+d-1}{d}}=\left\langle\begin{matrix}n\\k\end{matrix}\right\rangle_d$$

and

$$\prod_{j=0}^{d-1}\frac{\binom{n+d-1}{k+j}}{\binom{n+d-1}{j}}=\prod_{j=0}^{d-1}\frac{j!(n+d-1-j)!}{(k+j)!(n-k+d-1-j)!}=\prod_{j=0}^{d-1}\frac{j!(n+j)!}{(k+j)!(n-k+j)!}=\prod_{j=0}^{d-1}\frac{\binom{n+j}{k}}{\binom{k+j}{k}}.$$

An analog of (4) is

**Theorem 2**

$$\left\langle\begin{matrix}n\\k\end{matrix}\right\rangle_d=\det\left(\binom{n+i+j}{k+j}\right)_{i,j=0}^{d-1}=\det\left(\binom{n}{k+i-j}\right)_{i,j=0}^{d-1}. \qquad (8)$$

These determinants show that all $\left\langle\begin{matrix}n\\k\end{matrix}\right\rangle_d$ are integers.

**Proof**

Let us first prove the left-hand side.

$$\binom{n+i+j}{k+j}=\frac{(n+i+j)!}{(k+j)!(n+i-k)!}=\frac{j!(n+i)!}{(k+j)!(n+i-k)!}\frac{(n+i+j)!}{(j)!(n+i)!}=\frac{j!(n+i)!}{(k+j)!(n+i-k)!}\binom{n+i+j}{j}$$

implies $\det\left(\binom{n+i+j}{k+j}\right)_{i,j=0}^{d-1}=\prod_{j=0}^{d-1}\frac{j!(n+j)!}{(k+j)!(n+j-k)!}\det\left(\binom{n+i+j}{j}\right)_{i,j=0}^{d-1}$

with $\prod_{j=0}^{d-1}\frac{j!(n+j)!}{(k+j)!(n+j-k)!}=\prod_{j=0}^{d-1}\frac{\binom{n+j}{k}}{\binom{k+j}{k}}=\left\langle\begin{matrix}n\\k\end{matrix}\right\rangle_d$.

It remains to prove that

$$\det\left(\binom{n+i+j}{j}\right)_{i,j=0}^{d-1}=1. \tag{9}$$

To this end let $\Delta$ be the difference operator on the polynomials defined by $\Delta f(x)=f(x+1)-f(x)$. It satisfies $\Delta\binom{x}{n}=\binom{x+1}{n}-\binom{x}{n}=\binom{x}{n-1}$ and therefore $\Delta^k\binom{x}{n}=\binom{x}{n-k}$. Writing $\Delta=E-1$ with $Ef(x)=f(x+1)$ we get $\Delta^k=\left(E-1\right)^k=\sum_{j=0}^{k}(-1)^j\binom{k}{j}E^{k-j}$ and thus $\sum_{j=0}^{k}(-1)^j\binom{k}{j}\binom{x+k-j}{n}=\binom{x}{n-k}$.

Since $\binom{x}{n-k}=0$ for $k>n$ and $\binom{x}{n-n}=1$ the matrix $\left((-1)^{i-j}\binom{i}{j}\right)_{i,j=0}^{d-1}\left(\binom{n+i+j}{j}\right)_{i,j=0}^{d-1}$ is upper triangular with all entries $1$ in the main diagonal. This implies

$$\det\left(\binom{n+i+j}{j}\right)_{i,j=0}^{d-1}=\det\left((-1)^{i-j}\binom{i}{j}\right)_{i,j=0}^{d-1}\left(\binom{n+i+j}{j}\right)_{i,j=0}^{d-1}=1.$$

To compute the determinant $\det\left(\binom{n}{k+i-j}\right)_{i,j=0}^{d-1}=\det\left(\binom{n}{k+j-i}\right)_{i,j=1}^{d}$

we use formula (3.12) in [6] for $q=1$:

$$\det\left(\binom{A}{L_i+j}\right)_{i,j=1}^{d}=\frac{\prod_{1\le i<j\le n}\left(L_i-L_j\right)\prod_{i=1}^{n}(A+i-1)!}{\prod_{i=1}^{n}\left(L_i+n\right)!\prod_{i=1}^{n}\left(A-L_i-1\right)!}.$$

Choosing $A=n,\ L_i=k-i$ and $n=d$ this gives

$$\det\left(\binom{n}{k+j-i}\right)_{i,j=1}^{d}=\frac{\prod_{j=0}^{d-1}j!\prod_{j=0}^{d-1}(n+j)!}{\prod_{j=0}^{d-1}(k+j)!\prod_{j=0}^{d-1}(n-k+j)!}=\prod_{j=0}^{d-1}\frac{\binom{n+j}{k}}{\binom{k+j}{k}}=\left\langle{n \atop k}\right\rangle_d. \tag{10}$$

Another determinant representation has been given in [7]:

**Corollary 3**

$$\left\langle{n \atop k}\right\rangle_d=\det\left(\binom{n+i}{k+j}\right)_{i,j=0}^{d-1}. \tag{11}$$

**Proof**

If we subtract row $i-1$ from row $i$ in $\left(\binom{n+i}{k+j}\right)_{i,j=0}^{d-1}$ the new row $i$ has the entries

$\binom{n+i-1}{k+j-1}$. If we do this for $i=d-1,d-2,\cdots,1$ the new matrix has the first row unchanged

and the rest is the matrix $\binom{n+i-1}{k+j-1}_{i=1}^{d-1}$. If we iterate this we arrive at $\left(\binom{n}{k+j-i}\right)_{i,j=0}^{d-1}$.

**Remark**

In [14] and [7] these "MacMahon determinants" have been proved with the condensation method (cf. [6], Proposition 10). We will use this method in Theorem 8 for the proof of a $q-$ analog.

There is a nice generalization of the formula

$$\sum_{n\geq 0}\binom{n+k}{k}x^n=\frac{1}{(1-x)^{k+1}}. \tag{12}$$

In [12] Robert A. Sulanke introduced Narayana numbers $N(d,n,k)$ of dimension $d$. His results imply Theorem 4 which we state without proof.

**Theorem 4**

$$(1-x)^{dk+1}\sum_{n\geq 0}\left\langle {n+k \atop k}\right\rangle_d x^n=\sum_{j=0}^{(d-1)(k-1)}N(d,k,j)x^j. \tag{13}$$

For $d=3$ the polynomials $\sum_{j=0}^{2(k-1)}N(3,k,j)x^j$ are 1, $1+3x+x^2$, $1+10x+20x^2+10x^3+x^4$, $1+22x+113x^2+119x^3+113x^4+22x^5+x^6,\cdots$.

For $d=2$ we get the Narayana numbers $N(2,n,k)=\frac{1}{k+1}\binom{n}{k}\binom{n-1}{k}=N_{n,k}$ in the usual notation. In our notation $N_{n,k}=\left\langle {n-1 \atop k}\right\rangle_2$.

Let us give a direct proof for this case.

**Theorem 5**

$$\frac{\sum_{j=0}^{k-1}\left\langle {k-1 \atop j}\right\rangle_2 x^j}{(1-x)^{2k+1}}=\sum_{n\geq 0}\left\langle {n+k \atop k}\right\rangle_2 x^n. \tag{14}$$

**Proof**

Since $\binom{k+1}{j+1}=\frac{k+1}{j+1}\binom{k}{j}$ (14) is equivalent with

$$(1-x)^{2k+1}\sum_{n\geq 0}\binom{n+k}{k}\binom{n+k+1}{k}x^n=\sum_{j=0}^{k-1}\binom{k-1}{j}\binom{k+1}{j+1}x^j. \qquad (15)$$

If $D=\frac{d}{dx}$ denotes the differentiation operator we get

$$\sum_{n\geq 0}\binom{n+k}{k}\binom{n+k+1}{k}x^{n+1}=\frac{D^k}{k!}\sum_{n\geq 0}\binom{n+k}{k}x^{n+k+1}=\frac{D^k}{k!}\frac{x^{k+1}}{(1-x)^{k+1}}=\frac{D^k}{k!}\frac{\left(1-(1-x)\right)^{k+1}}{(1-x)^{k+1}}$$
$$=\frac{D^k}{k!}\sum_{j=0}^{k+1}(-1)^j\binom{k+1}{j}(1-x)^{j-k-1}=\sum_{j=0}^{k+1}(-1)^{j+k}\binom{k+1}{j}\binom{j-k-1}{k}(1-x)^{j-2k-1}$$
$$=\sum_{j=0}^{k+1}(-1)^j\binom{k+1}{j}\binom{2k-j}{k}(1-x)^{j-2k-1}.$$

It remains to show that

$$\sum_{j=0}^{k+1}(-1)^j\binom{k+1}{j}\binom{2k-j}{k}(1-x)^j=\sum_{j=0}^{k-1}\binom{k-1}{j}\binom{k+1}{j+1}x^{j+1}. \qquad (16)$$

Comparing the coefficient of $z^k$ in

$$\sum_{j,\ell}\binom{k-1}{j}\binom{k+1}{\ell}x^\ell z^{j+k+1-\ell}=(1+z)^{k-1}(x+z)^{k+1}=(1+z)^{k-1}\left(x-1+1+z\right)^{k+1}=\sum\binom{k+1}{j}(x-1)^j(1+z)^{2k-j}$$
$$=\sum_{j,\ell}\binom{k+1}{j}(x-1)^j\binom{2k-j}{\ell}z^\ell$$

gives (16) and thus (15).

## 3. A combinatorial interpretation

There is an interesting combinatorial interpretation which I owe to Qiaochu Yuan [13]:

$\left\langle{n \atop k}\right\rangle_d$ is the number of semistandard Young tableaux with shape $d^k$ (a box with $d$ columns and $k$ rows) and entries in $\{1,\cdots,n\}$. This is equivalent with all $k\times d-$matrices $\left(a_{i,j}\right)$ with entries in $\{1,\cdots,n\}$, such that $a_{i,j}\leq a_{i,j+1}$ and $a_{i,j}<a_{i+1,j}$ for all $i,j$.

For $d=1$ this is equivalent with choosing $k$ different numbers from $\{1,\cdots,n\}$.

In general the number of such matrices is given by the semistandard hook length formula (cf. [11]) which gives

$$\prod_{i=1}^{k}\prod_{j=1}^{d}\frac{n-i+j}{(k-i)+(d-j)+1}=\prod_{i=1}^{k}\frac{(n-i+1)(n-i+2)\cdots(n-i+d)}{(k-i+1)(k-i+2)\cdots(k-i+d)}$$
$$=\prod_{i=0}^{k-1}\frac{(n-i)(n-i+1)\cdots(n-i+d-1)}{(k-i)(k-i+1)\cdots(k-i+d-1)}=\prod_{i=0}^{k-1}\frac{\binom{n-i+d-1}{d}}{\binom{k-i+d-1}{d}}=\prod_{i=0}^{k-1}\frac{\langle n-i\rangle_d}{\langle k-i\rangle_d}=\left\langle {n \atop k}\right\rangle_d .$$

The Jacobi-Trudi identities (cf. [11]) give

**Theorem 6**

$$\left\langle {n \atop k}\right\rangle_d=\det\left(\binom{n+d+j-i-1}{n-1}\right)_{i,j=0}^{k-1}. \tag{17}$$

We now give an elementary

**Proof**

We consider more generally $\det\left(\binom{x_i+j}{n-1}\right)_{i,j=0}^{k-1}$. This is a polynomial in the indeterminates $x_0,\cdots,x_{k-1}$ of degree $\le n-1$ in each variable. It vanishes for $x_i=x_j$ which gives the factor $\prod_{0\le i<j\le k-1}(x_j-x_i).$ Since all entries of row $i$ have the factor $x_i(x_i-1)\cdots(x_i-n+k+1)$ the determinant has the factor $\prod_{0\le i<j\le k-1}(x_j-x_i)\prod_{i=0}^{k-1}x_i(x_i-1)\cdots(x_i-n+k+1).$

This also is a polynomial of degree $n$ in each variable. Therefore, there exists a constant $c$ such that

$$\det\left(\binom{x_i+j}{n-1}\right)_{i,j=0}^{k-1}=c\prod_{0\le i<j\le k-1}(x_j-x_i)\prod_{i=0}^{k-1}\binom{x_i}{n-k}. \tag{18}$$

To compute $c$ we choose $x_i=n-1-i.$ Then $\left(\binom{x_i+j}{n-1}\right)_{i,j=0}^{k-1}$ is a right triangle matrix with $\binom{x_i+i}{n-1}=\binom{n-1}{n-1}=1$ und therefore $\det\left(\binom{x_i+j}{n-1}\right)_{i,j=0}^{k-1}=1.$

On the right-hand side of (18) we get

$$c\prod_{0\le i<j\le k-1}(x_j-x_i)\prod_{i=0}^{k-1}\binom{x_i}{n-k}=c\prod_{0\le i<j\le k-1}(i-j)\prod_{i=0}^{k-1}\binom{x_i}{n-k}$$
$$=c(-1)^{\binom{k}{2}}\prod_{i=0}^{k-1}i!\prod_{i=0}^{k-1}\binom{n-1-i}{k-1-i}$$

Setting $f(k)=\prod_{i=0}^{k-1} i!\binom{n-1-i}{k-1-i}$ we get

$$\frac{f(k)}{f(k-1)}=\frac{\prod_{i=0}^{k-1} i!\binom{n-1-i}{k-1-i}}{\prod_{i=0}^{k-2} i!\binom{n-1-i}{k-2-i}}=(k-1)!\prod_{i=0}^{k-2}\frac{(k-2-i)!(n-k+1)!}{(k-1-i)!(n-k)!}=$$

$$(k-1)!\prod_{i=0}^{k-2}\frac{n-k+1}{k-1-i}=(n-k+1)^k$$

and therefore $f(k)=\prod_{j=0}^{k-1}(n-j)^j$.

Thus

$$\det\left(\binom{x_i+j}{n-1}\right)_{i,j=0}^{k-1}=(-1)^{\binom{k}{2}}\frac{1}{\prod_{j=0}^{k-1}(n-j)^j}\prod_{0\le i<j\le k-1}(x_j-x_i)\prod_{i=0}^{k-1}\binom{x_i}{n-k}. \tag{19}$$

To compute $\det\left(\binom{d-i+n-1+j}{n-1}\right)_{i,j=1}^{k}$ we choose $x_i=d-i+n-1$ and get

$$(-1)^{\binom{k}{2}}\frac{1}{\prod_{j=0}^{k-1}(n-j)^j}\prod_{0\le i<j\le k-1}(i-j)\prod_{i=0}^{k-1}\binom{d-i+n-1}{n-k}=\prod_{j=0}^{k-1}\frac{j!}{(n-j)^j}\frac{(d-i+n-1)!}{(n-k)!(d-i+k-1)!}$$

$$=\prod_{i=0}^{k-1}\frac{\binom{n-i+d-1}{d}}{\binom{k-i+d-1}{d}}\prod_{i=0}^{k-1}\frac{(n-i-1)!i!}{(k-i-1)!(n-i)^i(n-k)!}=\left\langle {n \atop k} \right\rangle_d,$$

because

$$\prod_{i=0}^{k-1}\frac{(n-i-1)!i!}{(k-i-1)!(n-i)^i(n-k)!}=\prod_{i=0}^{k-1}\frac{(n-i-1)!}{(n-i)^i(n-k)!}$$

$$=\frac{(n-1)!(n-2)!\cdots(n-k)!}{(n-1)(n-2)^2\cdots(n-k+1)^{k-1}(n-k)^k(n-k-1)!^k}=1.$$

## 4. q-analogs

The above constructions have straightforward $q-$analogs. For a real number $q$ with $|q|<1$

let $[n]_q=1+q+\cdots+q^{n-1}=\frac{1-q^n}{1-q}$, $[n]_q!=\prod_{j=1}^{n}[j]_q$ and $\begin{bmatrix} n \\ k \end{bmatrix}_q=\frac{[n]_q!}{[k]_q![n-k]_q!}=\prod_{j=0}^{k-1}\frac{1-q^{n-j}}{1-q^{k-j}}$.

As is well known the $q-$binomial coefficients $\begin{bmatrix} n \\ k \end{bmatrix}_q$ satisfy

$$\begin{bmatrix} n \\ k \end{bmatrix}_q = q^k \begin{bmatrix} n-1 \\ k \end{bmatrix}_q + \begin{bmatrix} n-1 \\ k-1 \end{bmatrix}_q = \begin{bmatrix} n-1 \\ k \end{bmatrix}_q + q^{n-k} \begin{bmatrix} n-1 \\ k-1 \end{bmatrix}_q \tag{20}$$

and are therefore polynomials in $q$ with nonnegative integer coefficients.

For later use let us mention the following $q-$analog of the binomial theorem (cf. e.g. [1])

$$\prod_{j=0}^{n-1}(1-q^j x) = \sum_{j=0}^{n}(-1)^j q^{\binom{j}{2}} \begin{bmatrix} n \\ j \end{bmatrix}_q x^j. \tag{21}$$

We define $\langle n \rangle_{d,q} = \begin{bmatrix} n+d-1 \\ d \end{bmatrix}_q$, $\langle n \rangle_{d,q}! = \prod_{j=1}^{n} \langle j \rangle_{d,q}$ and get

$$\left\langle \begin{matrix} n \\ k \end{matrix} \right\rangle_{d,q} = \frac{\langle n \rangle_{d,q}!}{\langle k \rangle_{d,q}! \langle n-k \rangle_{d,q}!} = \prod_{j=0}^{k-1} \frac{\langle n-j \rangle_{d,q}}{\langle k-j \rangle_{d,q}} = \prod_{j=0}^{d-1} \frac{\begin{bmatrix} n+j \\ k \end{bmatrix}_q}{\begin{bmatrix} k+j \\ k \end{bmatrix}_q} = \prod_{j=0}^{d-1} \frac{\begin{bmatrix} n+d-1 \\ k+j \end{bmatrix}_q}{\begin{bmatrix} n+d-1 \\ j \end{bmatrix}_q}, \tag{22}$$

From Theorem 7 we see that these are also polynomials in $q$ with integer coefficients.

For $d=2$ we get $\left\langle \begin{matrix} n \\ k \end{matrix} \right\rangle_{2,q} = \frac{1}{[k+1]_q} \begin{bmatrix} n \\ k \end{bmatrix}_q \begin{bmatrix} n+1 \\ k \end{bmatrix}_q$.

This gives the triangle

$$\begin{pmatrix} 1 & & & & \\ 1 & 1 & & & \\ 1 & 1+q+q^2 & 1 & & \\ 1 & (1+q^2)(1+q+q^2) & (1+q^2)(1+q+q^2) & 1 & \\ 1 & (1+q^2)(1+q+q^2+q^3+q^4) & (1+q^2)^2(1+q+q^2+q^3+q^4) & (1+q^2)(1+q+q^2+q^3+q^4) & 1 \end{pmatrix}$$

As analog of (8) we get

**Theorem 7**

$$\left\langle \begin{matrix} n \\ k \end{matrix} \right\rangle_{d,q} = \frac{\det\left( \begin{bmatrix} n+i+j \\ k+j \end{bmatrix}_q \right)_{i,j=0}^{d-1}}{q^{n\binom{d}{2} + \sum_{j=0}^{d-1} j^2}} = \det\left( q^{\binom{i-j}{2}} \begin{bmatrix} n \\ k-i+j \end{bmatrix}_q \right)_{i,j=0}^{d-1}. \tag{23}$$

**Proof**

From $\begin{bmatrix} n+i+j \\ k+j \end{bmatrix}_q = \frac{[n+i+j]_q!}{[k+j]_q![n-k+i]_q!} = \frac{[j]_q![n+i]_q!}{[k+j]_q![n-k+i]_q!} \begin{bmatrix} n+i+j \\ j \end{bmatrix}_q$ we get

$$\det\left(\begin{bmatrix} n+i+j \\ k+j \end{bmatrix}_q\right)_{i,j=0}^{d-1} = \prod_{j=0}^{d-1} \frac{[j]_q![n+j]_q!}{[k+j]_q![n-k+j]_q!} \det\left(\begin{bmatrix} n+i+j \\ j \end{bmatrix}_q\right)_{i,j=0}^{d-1}$$

with $\prod_{j=0}^{d-1} \frac{[j]_q![n+j]_q!}{[k+j]_q![n-k+j]_q!} = \prod_{j=0}^{d-1} \frac{\begin{bmatrix} n+j \\ k \end{bmatrix}_q}{\begin{bmatrix} k+j \\ k \end{bmatrix}_q} = \left\langle \begin{matrix} n \\ k \end{matrix} \right\rangle_{d,q}$.

It remains to compute $\det\left(\begin{bmatrix} n+i+j \\ j \end{bmatrix}_q\right)_{i,j=0}^{d-1}$.

Using the identity

$$\sum_{\ell=0}^{i} (-1)^{i-\ell} \begin{bmatrix} i \\ \ell \end{bmatrix}_q q^{\binom{i-\ell}{2}} \begin{bmatrix} n+\ell+j \\ j \end{bmatrix}_q = q^{i(n+i)} \begin{bmatrix} n+j \\ j-i \end{bmatrix}_q \tag{24}$$

we see that $\left((-1)^{i-\ell} q^{\binom{i-\ell}{2}} \begin{bmatrix} i \\ \ell \end{bmatrix}\right)_{i,\ell=0}^{d-1} \left(\begin{bmatrix} n+i+j \\ j \end{bmatrix}_q\right)_{i,j=0}^{d-1}$ is an upper triangular matrix with entries $q^{in+i^2}$ in the main diagonal. This implies

$$\det\left(\begin{bmatrix} n+i+j \\ j \end{bmatrix}_q\right)_{i,j=0}^{d-1} = q^{n\binom{d}{2}+\sum_{j=0}^{d-1} j^2}. \tag{25}$$

To prove (24) we consider $\begin{bmatrix} x \\ k \end{bmatrix}_q = \prod_{j=0}^{k-1} \frac{[x-j]_q}{[k-j]_q} = \prod_{j=0}^{k-1} \frac{q^{x-j}-1}{q^{k-j}-1} = \prod_{j=0}^{k-1} \frac{q^x - q^j}{q^k - q^j}$ as a polynomial in $q^x$ with coefficients in $Q(q)$. If we define the operator $E$ on these polynomials by $Ef(q^x) = f\left(q^{x+1}\right)$ then we get $E\begin{bmatrix} x \\ n \end{bmatrix} = \begin{bmatrix} x+1 \\ n \end{bmatrix}$. Let now $\Delta = E - 1$ be the difference operator.

We have $\Delta\begin{bmatrix} x \\ n \end{bmatrix} = \begin{bmatrix} x+1 \\ n \end{bmatrix} - \begin{bmatrix} x \\ n \end{bmatrix} = q^{x-n+1}\begin{bmatrix} x \\ n-1 \end{bmatrix}$.

More generally we get by induction

$$(E-1)(E-q)\cdots\left(E-q^{k-1}\right)\begin{bmatrix} x \\ n \end{bmatrix}_q = q^{k(x+k-n)}\begin{bmatrix} x \\ n-k \end{bmatrix}_q,$$

because

$$\left(E-q^{k-1}\right)q^{(k-1)(x+k-1-n)}\begin{bmatrix} x \\ n-k+1 \end{bmatrix}_q = q^{(k-1)(x+k-n)}\begin{bmatrix} x+1 \\ n-k+1 \end{bmatrix}_q - q^{(k-1)(x+k-n)}\begin{bmatrix} x \\ n-k+1 \end{bmatrix}_q$$

$$= q^{(k-1)(x+k-n)}q^{x+k-n}\begin{bmatrix} x \\ n-k \end{bmatrix}_q = q^{k(x+k-n)}\begin{bmatrix} x \\ n-k \end{bmatrix}_q .$$

Finally by (21) we have $(E-1)(E-q)\cdots\left(E-q^{k-1}\right)=\sum_{j=0}^{k}(-1)^j q^{\binom{j}{2}}\begin{bmatrix} k \\ j \end{bmatrix}E^{k-j}.$

For the computation of $\det\left(q^{\binom{i-j}{2}}\begin{bmatrix} n \\ k-j+i \end{bmatrix}_q\right)_{i,j=0}^{d-1} = \det\left(q^{\binom{j-i}{2}}\begin{bmatrix} n \\ k-i+j \end{bmatrix}_q\right)_{i,j=1}^{d}$ we use [6], formula (3.12):

$$\det\left(q^{jL_i}\begin{bmatrix} A \\ L_i+j \end{bmatrix}_q\right)_{i,j=1}^{n} = q^{\sum_{i=1}^{n} iL_i}\frac{\prod_{1\le i<j\le n}\left[L_i-L_j\right]_q \prod_{i=1}^{n}[A+i-1]_q!}{\prod_{i=1}^{n}[L_i+n]_q!\prod_{i=1}^{n}[A-L_i-1]_q!}.$$

First we write

$$\det\left(q^{\binom{j-i}{2}}\begin{bmatrix} n \\ k-i+j \end{bmatrix}_q\right)_{i,j=1}^{d} = \det\left(q^{\binom{i+1}{2}+\binom{j}{2}-ij}\begin{bmatrix} n \\ k-i+j \end{bmatrix}_q\right)_{i,j=1}^{d} = q^{\sum_{j=1}^{d} j^2}\det\left(q^{-ij}\begin{bmatrix} n \\ k-i+j \end{bmatrix}_q\right)_{i,j=1}^{d}$$

$$\det\left(q^{\binom{j-i}{2}}\begin{bmatrix} n \\ k-i+j \end{bmatrix}_q\right)_{i,j=1}^{d} = q^{\sum_{j=1}^{d} j^2-k\binom{d+1}{2}}\det\left(q^{(k-i)j}\begin{bmatrix} n \\ k-i+j \end{bmatrix}_q\right)_{i,j=1}^{d}.$$

Then we choose $A=n,\ L_i=k-i$, $n=d$ and get

$$q^{\sum_{j=1}^{d} j^2-k\binom{d+1}{2}}\det\left(q^{(k-i)j}\begin{bmatrix} n \\ k-i+j \end{bmatrix}_q\right)_{i,j=1}^{d} = q^{\sum_{j=0}^{d} j^2-k\binom{d+1}{2}}q^{\sum_{i=1}^{d}(k-i)i}\frac{\prod_{j=0}^{d-1}[j]_q!\prod_{j=0}^{d-1}[n+j]_q!}{\prod_{j=0}^{d-1}[k+j]_q!\prod_{j=0}^{d-1}[n-k+j]_q!}$$

$$=\frac{\prod_{j=0}^{d-1}[j]_q!\prod_{j=0}^{d-1}[n+j]_q!}{\prod_{j=0}^{d-1}[k+j]_q!\prod_{j=0}^{d-1}[n-k+j]_q!}=\prod_{j=0}^{d-1}\frac{\begin{bmatrix} n+j \\ k \end{bmatrix}_q}{\begin{bmatrix} k+j \\ k \end{bmatrix}_q}=\left\langle \begin{matrix} n \\ k \end{matrix}\right\rangle_{d,q}.$$

**Theorem 8**

$$\left\langle \begin{matrix} n \\ k \end{matrix} \right\rangle_{d,q} = q^{-(n-k)\binom{d}{2}} \det\left( \begin{bmatrix} n+i \\ k+j \end{bmatrix}_q \right)_{i,j=0}^{d-1}. \tag{26}$$

**Proof**

We use Dodgson's condensation method (cf. [6], Proposition 10, and [14]). Let $X(d,n,k) = \det\left( \begin{bmatrix} n+i \\ k+j \end{bmatrix}_q \right)_{i,j=0}^{d-1}$. By condensation we get

$$X(d,n,k) = \frac{X(d-1,n,k)X(d-1,n+1,k+1) - X(d-1,n+1,k)X(d-1,n,k+1]}{X(d-2,n+1,k+1)}.$$

The same identity holds for $X(d,n,k) = q^{(n-k)\binom{d}{2}} \left\langle \begin{matrix} n \\ k \end{matrix} \right\rangle_{d,q}$. For we have

$$\frac{q^{(n-k)\binom{d-1}{2}} \left\langle \begin{matrix} n \\ k \end{matrix} \right\rangle_{d-1,q} q^{(n-k)\binom{d-1}{2}} \left\langle \begin{matrix} n+1 \\ k+1 \end{matrix} \right\rangle_{d-1,q}}{q^{(n-k)\binom{d}{2}} \left\langle \begin{matrix} n \\ k \end{matrix} \right\rangle_{d,q} q^{(n-k)\binom{d-2}{2}} \left\langle \begin{matrix} n+1 \\ k+1 \end{matrix} \right\rangle_{d-2,q}} = \frac{1}{q^{n-k}} \frac{[n-k+d-1]_q}{[d-1]_q}$$

and

$$\frac{q^{(n-k)\binom{d-1}{2}} \left\langle \begin{matrix} n+1 \\ k \end{matrix} \right\rangle_{d-1,q} q^{(n-k)\binom{d-1}{2}} \left\langle \begin{matrix} n \\ k+1 \end{matrix} \right\rangle_{d-1,q}}{q^{(n-k)\binom{d}{2}} \left\langle \begin{matrix} n \\ k \end{matrix} \right\rangle_{d,q} q^{(n-k)\binom{d-2}{2}} \left\langle \begin{matrix} n+1 \\ k+1 \end{matrix} \right\rangle_{d-2,q}} = \frac{1}{q^{n-k}} \frac{[n-k]_q}{[d-1]_q},$$

which implies

$$\frac{X(d-1,n,k)X(d-1,n+1,k+1) - X(d-1,n+1,k)X(d-1,n,k+1]}{X(d,n,k)X(d-2,n+1,k+1)} = \frac{[n-k+d-1]_q - [n-k]_q}{q^{n-k}[d-1]_q} = 1.$$

Identity (26) holds for $d=0$ and $d=1$ and therefore (26) holds for all $d$ by induction.

A $q$-analog of (17) is

**Theorem 9**

$$\det\left( \begin{bmatrix} n+d+j-i-1 \\ n-1 \end{bmatrix}_q \right)_{i,j=0}^{k-1} = q^{\binom{k}{2}d} \left\langle \begin{matrix} n \\ k \end{matrix} \right\rangle_{d,q}. \tag{27}$$

**Proof**

Here we use [6], formula (3.11):

$$\det\left(\begin{bmatrix} L_i+A+j \\ L_i+j \end{bmatrix}_q\right)_{i,j=1}^{n} = q^{\sum_{i=1}^{n}(i-1)(L_i+i)} \frac{\prod_{1\le i<j\le n}\left[L_i-L_j\right]_q \prod_{i=1}^{n}\left[L_i+A+1\right]_q!}{\prod_{i=1}^{n}\left[L_i+n\right]_q! \prod_{i=1}^{n}\left[A+1-i\right]_q!}.$$

We choose $L_i = d-i,\ A = n-1$ and $n = k$ and get

$$\det\left(\begin{bmatrix} n+d+j-i-1 \\ n-1 \end{bmatrix}_q\right)_{i,j=0}^{k-1} = q^{\binom{k}{2}d} \prod_{j=1}^{k} \frac{\left[j-1\right]_q!\left[d+n-j\right]_q!}{\left[d+k-j\right]_q!\left[n-j\right]_q!} = q^{\binom{k}{2}d} \prod_{j=0}^{k-1} \frac{\left[k-1-j\right]_q!\left[d+n-1-j\right]_q!}{\left[d+k-1-j\right]_q!\left[n-1-j\right]_q!}$$

$$= q^{\binom{k}{2}d} \prod_{j=0}^{k-1} \frac{\begin{bmatrix} n-j+d-1 \\ d \end{bmatrix}_q}{\begin{bmatrix} k-j+d-1 \\ d \end{bmatrix}_q} = q^{\binom{k}{2}d} \prod_{j=0}^{d-1} \frac{\begin{bmatrix} n+j \\ k \end{bmatrix}_q}{\begin{bmatrix} k+j \\ k \end{bmatrix}_q} = q^{\binom{k}{2}d} \left\langle \begin{matrix} n \\ k \end{matrix} \right\rangle_{d,q}.$$

**Remark**

Since the determinant (27) is closely related to semistandard Young Tableaux it would make sense from this point of view to define $q-$Hoggatt matrices with entries $q^{\binom{k}{2}d} \left\langle \begin{matrix} n \\ k \end{matrix} \right\rangle_{d,q}$ instead of $\left\langle \begin{matrix} n \\ k \end{matrix} \right\rangle_{d,q}$. For $d=1$ this means to replace the $q-$binomial coefficients $\begin{bmatrix} n \\ k \end{bmatrix}_q$ by their companion form $q^{\binom{k}{2}} \begin{bmatrix} n \\ k \end{bmatrix}_q$. This would give the nice generating function

$\sum_{k=0}^{n} q^{\binom{k}{2}} \begin{bmatrix} n \\ k \end{bmatrix}_q x^k = (1+x)(1+qx)\cdots\left(1+q^{n-1}x\right)$ for row $n$ of the matrix.

For $d=2$ we would get

$$\sum_{k=0}^{n} q^{k(k+1)} \left\langle \begin{matrix} n \\ k \end{matrix} \right\rangle_{2,q} = C_{n+1}(q) = \frac{1}{[n+2]} \begin{bmatrix} 2n+2 \\ n+1 \end{bmatrix},$$

which is a nice $q-$analog of the fact that the sum of the Narayana numbers $\left\langle \begin{matrix} n \\ k \end{matrix} \right\rangle_2$ are the Catalan numbers $C_{n+1}$.

As $q-$analog of (13) we state

**Conjecture 10**

*For positive integers* $d,k$ *we have*

$$(1-x)(1-qx)\cdots\left(1-q^{dk}x\right)\sum_{n\geq 0}\left\langle {n+k \atop k} \right\rangle_{d,q} x^n = \sum_{j=0}^{(d-1)(k-1)} N(d,k,j,q)x^j. \quad (28)$$

*where the coefficients* $N(d,k,j,q)$ *are palindromic polynomials in* $q$ *with nonnegative coefficients.*

*For* $d=2$ *this reduces to*

$$(1-x)(1-qx)\cdots\left(1-q^{2k}x\right)\sum_{n\geq 0}\left\langle {n+k \atop k} \right\rangle_{2,q} x^n = \sum_{j=0}^{k-1} q^{j(j+1)}\left\langle {k-1 \atop j} \right\rangle_{2,q} x^j. \quad (29)$$

*The sums* $\sum_{j=0}^{(d-1)(k-1)} N(d,k,j,q) = C_n^{(d)}(q) = [dn]_q!\prod_{j=0}^{d-1}\frac{[j]_q!}{[n+j]_q!}$ are the $d-$dimensional $q-$Catalan numbers.

## 5. Fibonacci-Hoggatt triangles

Let $F_n = \sum_{j=0}^{\left\lfloor \frac{n-1}{2} \right\rfloor}\binom{n-1-j}{j}$ denote the Fibonacci numbers which satisfy $F_n = F_{n-1} + F_{n-2}$ with initial values $F_0 = 0$ and $F_1 = 1$ and $F_n = \frac{\alpha^n - \beta^n}{\alpha - \beta}$ with $\alpha = \frac{1+\sqrt{5}}{2}$ and $\beta = \frac{1-\sqrt{5}}{2}$.

Let us write $(n)_F = F_n$, $(n)_F! = F_nF_{n-1}\cdots F_1$ and define the Fibonomial coefficients by $\binom{n}{k}_F = \frac{F_nF_{n-1}\cdots F_{n-k+1}}{F_kF_{k-1}\cdots F_1} = \frac{(n)_F!}{(k)_F!(n-k)_F!}$.

The first terms are (cf. OEIS [8], A010048)

$$\begin{pmatrix} 1 & & & & & \\ 1 & 1 & & & & \\ 1 & 1 & 1 & & & \\ 1 & 2 & 2 & 1 & & \\ 1 & 3 & 6 & 3 & 1 & \\ 1 & 5 & 15 & 15 & 5 & 1 \end{pmatrix}.$$

The Fibonacci numbers satisfy

$$\begin{pmatrix} 0 & 1 \\ 1 & 1 \end{pmatrix}^n = \begin{pmatrix} F_{n-1} & F_n \\ F_n & F_{n+1} \end{pmatrix}.$$

From $\begin{pmatrix} F_{n-1} & F_n \\ F_n & F_{n+1} \end{pmatrix} = \begin{pmatrix} 0 & 1 \\ 1 & 1 \end{pmatrix}^{n-k} \begin{pmatrix} 0 & 1 \\ 1 & 1 \end{pmatrix}^{k} = \begin{pmatrix} F_{n-k-1} & F_{n-k} \\ F_{n-k} & F_{n-k+1} \end{pmatrix} \begin{pmatrix} F_{k-1} & F_k \\ F_k & F_{k+1} \end{pmatrix}$

we get by comparing the top right elements

$F_n = F_{k+1}F_{n-k} + F_{n-k-1}F_k.$

This is equivalent with

$$\binom{n}{k}_F = F_{k+1}\binom{n-1}{k}_F + F_{n-k-1}\binom{n-1}{k-1}_F, \tag{30}$$

which shows that the Fibonomials $\binom{n}{k}_F$ are nonnegative integers.

L. Carlitz [3] found the analog of (21)

$$h_n(x) = \sum_{j=0}^{n} (-1)^{\binom{j+1}{2}} \binom{n}{j}_F x^j = \prod_{j=0}^{n-1} \left(1 - \alpha^{n-j-1}\beta^j x\right). \tag{31}$$

Let us reproduce his proof. In formula (21) we set $q = \frac{\beta}{\alpha}$ and get

$$\prod_{j=0}^{n-1}\left(1 - \frac{\beta^j}{\alpha^j}x\right) = \sum_{j=0}^{n} (-1)^j \left(\frac{\beta}{\alpha}\right)^{\binom{j}{2}} \begin{bmatrix} n \\ j \end{bmatrix}_{\frac{\beta}{\alpha}} x^j$$

with

$$\begin{bmatrix} n \\ j \end{bmatrix}_q = \frac{\left(1-q^n\right)\cdots\left(1-q^{n-j+1}\right)}{\left(1-q^j\right)\cdots\left(1-q\right)} = \frac{1}{\alpha^{(n-j)j}} \frac{\left(\alpha^n - \beta^n\right)\cdots\left(\alpha^{n-j+1} - \beta^{n-j+1}\right)}{\left(\alpha^j - \beta^j\right)\cdots\left(\alpha - \beta\right)} = \alpha^{j^2-nj}\binom{n}{j}_F.$$

This gives

$$\prod_{j=0}^{n-1}\left(1 - \alpha^{-j}\beta^j x\right) = \sum_{j=0}^{n} (-1)^j \alpha^{-\binom{j}{2}} \beta^{\binom{j}{2}} \alpha^{j^2-nj} \binom{n}{j}_F x^j.$$

Replacing $x \to \alpha^{n-1}x$ we get

$$\prod_{j=0}^{n-1}\left(1 - \alpha^{n-j-1}\beta^j x\right) = \sum_{j=0}^{n} (-1)^j (\alpha\beta)^{\binom{j}{2}} \binom{n}{j}_F x^j = \sum_{j=0}^{n} (-1)^{\binom{j+1}{2}} \binom{n}{j}_F x^j.$$

Since $\alpha^n + \beta^n = L_n$ are the Lucas numbers $2,1,3,4,7,11,18,\cdots,$ we see that

$$h_n(x) = \prod_{j=0}^{n-1}\left(1 - \alpha^{n-j-1}\beta^j x\right) = \left(1-\alpha^{n-1}x\right)\left(1-\beta^{n-1}x\right)\prod_{j=1}^{n-2}\left(1 + \alpha^{n-2-j}\beta^{j-1}x\right)$$

$$= \left(1-\alpha^{n-1}x\right)\left(1-\beta^{n-1}x\right)\prod_{j=0}^{n-3}\left(1 + \alpha^{n-3-j}\beta^{j}x\right) = \left(1 - L_{n-1}x + (-1)^{n-1}x^2\right)p_{n-2}(-x).$$

This gives

$$h_k(x)=\sum_{j=0}^{k}(-1)^{\binom{j+1}{2}}\binom{k}{j}_F x^j=\prod_{j=0}^{\left\lfloor \frac{k}{2}\right\rfloor}u_{k-2j}\left((-1)^j x\right) \tag{32}$$

with

$$\begin{aligned}&u_k(x)=1-L_{k-1}x+(-1)^{k-1}x^2,\\&u_1(x)=1-x,\quad u_0(x)=1.\end{aligned} \tag{33}$$

As analog of (12) we get

$$\frac{1}{h_{k+1}(x)}=\frac{1}{\sum_{j=0}^{k+1}(-1)^{\binom{j+1}{2}}\binom{k+1}{j}_F x^j}=\sum_{n\geq 0}\binom{n+k}{k}_F x^n. \tag{34}$$

**Proof**

Since $h_1(x)=1-x$ and $h_2(x)=1-x-x^2$ identity (34) is true for $k=0$ and $k=1$.

By (32) identity (34) is equivalent with $u_{k+1}(x)\sum_{n\geq 0}\binom{n+k}{k}_F x^n=\sum_{n\geq 0}\binom{n+k-2}{k-2}_F(-x)^n$, i.e.

$$\left(1-L_k x+(-1)^k x^2\right)\sum_{n\geq 0}\frac{F_{n+1}\cdots F_{n+k}}{F_1\cdots F_k}x^n=\sum_{n\geq 0}\frac{F_{n+1}\cdots F_{n+k-2}}{F_1\cdots F_{k-2}}(-x)^n.$$

This is equivalent with

$$F_{n+k-1}F_{n+k}-L_kF_nF_{n+k-1}+(-1)^kF_{n-1}F_n=(-1)^nF_{k-1}F_k,$$

which is easily verified.

In [5] the authors studied Fibo-Narayana numbers defined by

$$\left\langle {n \atop k} \right\rangle_{2,F}=\frac{1}{F_{k+1}}\binom{n}{k}_F\binom{n+1}{k}_F=\frac{1}{F_{n+1}}\binom{n+1}{k}_F\binom{n+1}{k+1}_F.$$

Let us more generally define Fibo-Hoggatt numbers

$$\left\langle {n \atop k} \right\rangle_{d,F}=\prod_{j=0}^{k-1}\frac{\langle n-j\rangle_{d,F}}{\langle k-j\rangle_{d,F}} \tag{35}$$

with $\langle n\rangle_{d,F}=\binom{n+d-1}{d}_F$ and consider the corresponding Fibonacci-Hoggatt matrices

$$H_{d,F}=\left(\left\langle {n \atop k}\right\rangle_{d,F}\right)_{n,k\geq 0}.$$

As in (7) we get

$$\left\langle {n \atop k}\right\rangle_{d,F}=\prod_{j=0}^{d-1}\frac{\binom{n+j}{k}_F}{\binom{k+j}{k}_F}=\prod_{j=0}^{d-1}\frac{\binom{n+d-1}{k+j}_F}{\binom{n+d-1}{j}_F}. \tag{36}$$

For example for $d=3$ we get

$$\left(\left\langle {n \atop k}\right\rangle_{3,F}\right)_{i.j=0}^{5}=\begin{pmatrix}1 & & & & & \\ 1 & 1 & & & & \\ 1 & 3 & 1 & & & \\ 1 & 15 & 15 & 1 & & \\ 1 & 60 & 300 & 60 & 1 & \\ 1 & 260 & 5200 & 5200 & 260 & 1\end{pmatrix}.$$

As an analog of the first identity (8) we get

**Theorem 11**

$$\left\langle {n \atop k}\right\rangle_{d,F}=\frac{\det\left(\binom{n+i+j}{k+j}_F\right)_{i,j=0}^{d-1}}{\det\left(\binom{n+i+j}{j}_F\right)_{i,j=0}^{d-1}}. \tag{37}$$

**Proof**

This follows from

$$\det\left(\binom{n+i+j}{k+j}_F\right)_{i,j=0}^{d-1}=\det\left(\frac{(j)_F!(n+i)_F!}{(k+j)_F!(n-k+i)_F!}\binom{n+i+j}{j}_F\right)_{i,j=0}^{d-1}$$

$$=\prod_{j=0}^{d-1}\frac{(j)_F!(n+j)_F!}{(k+j)_F!(n-k+j)_F!}\det\left(\binom{n+i+j}{j}_F\right)_{i,j=0}^{d-1}$$

if we observe that

$$\prod_{j=0}^{d-1}\frac{(j)_F!(n+j)_F!}{(k+j)_F!(n-k+j)_F!}=\prod_{j=0}^{d-1}\frac{\binom{n+j}{k}_F}{\binom{k+j}{k}_F}=\left\langle {n \atop k} \right\rangle_{d,F}.$$

If we set $a(d,n)=\det\left(\binom{n+i+j}{j}_F\right)_{i,j=0}^{d-1}$ then we get $a(2,n)=F_n$, but for $n>3$ no other interpretation seems to be known. For example for $d=3$ we get $1,5,7,53,187,853,\cdots$.

As analog of (10) we get as special case of Theorem 15

**Theorem 12**

$$\det\left((-1)^{\binom{i-j}{2}}\binom{n}{k-i+j}_F\right)_{i,j=0}^{d-1}=\left\langle {n \atop k} \right\rangle_{d,F}. \tag{38}$$

This shows that all $\left\langle {n \atop k} \right\rangle_{d,F}$ are positive integers.

There is also a nice analog of (17):

**Theorem 13**

$$\left\langle {n \atop k} \right\rangle_{d,F}=(-1)^{d\binom{k}{2}}\det\left(\binom{n+d+j-i-1}{n-1}_F\right)_{i,j=0}^{k-1}. \tag{39}$$

**Proof**

By Binet's formula we have $F_n=\alpha^{n-1}\frac{1-\left(\frac{\beta}{\alpha}\right)^n}{1-\left(\frac{\beta}{\alpha}\right)}$ with $\frac{\beta}{\alpha}=-\frac{3-\sqrt{5}}{2}=-\frac{1}{\alpha^2}$. If we set $\frac{\beta}{\alpha}=q$, then we get

$$F_n=(-q)^{\frac{1-n}{2}}\frac{1-q^n}{1-q}. \tag{40}$$

This implies

$$\binom{n}{k}_F=(-q)^{\frac{k^2-nk}{2}}\begin{bmatrix} n \\ k \end{bmatrix}_q. \tag{41}$$

We can now formulate $\det\left(\binom{n+d+j-i-1}{n-1}_F\right)_{i,j=0}^{k-1}$ in terms of $q-$binomial coefficients as

$$\det\left(\binom{n+d+j-i-1}{n-1}_F\right)_{i,j=0}^{k-1}=\det\left((-q)^{\frac{(1-n)(d+j-i)}{2}}\begin{bmatrix}n+d+j-i-1\\ n-1\end{bmatrix}\right)_{i,j=0}^{k-1}$$

$$=(-q)^{\frac{(1-n)dk}{2}}\det\left(\begin{bmatrix}n+d+j-i-1\\ d+j-i\end{bmatrix}\right)_{i,j=0}^{k-1}.$$

Using the above result we get by (41)

$$\det\left(\binom{n+d+j-i-1}{n-1}_F\right)_{i,j=0}^{k-1}=(-q)^{\frac{(1-n)dk}{2}}q^{\binom{k}{2}d}\prod_{j=0}^{k-1}\frac{\begin{bmatrix}n-j+d-1\\ d\end{bmatrix}_q}{\begin{bmatrix}k-j+d-1\\ d\end{bmatrix}_q}$$

$$=(-q)^{\frac{(1-n)dk}{2}}q^{\binom{k}{2}d}(-q)^{\frac{dnk-dk^2}{2}}\prod_{j=0}^{k-1}\frac{\binom{n-j+d-1}{d}_F}{\binom{k-j+d-1}{d}_F}=(-1)^{\binom{k}{2}d}\prod_{j=0}^{k-1}\frac{\binom{n-j+d-1}{d}_F}{\binom{k-j+d-1}{d}_F}=(-1)^{\binom{k}{2}d}\left\langle{n\atop k}\right\rangle_{d,F}$$

As analogs of (28) and (29) we state

**Conjecture 14**

$\left(\sum_{j=0}^{dk+1}(-1)^{\binom{j+1}{2}}\binom{dk+1}{j}_F x^j\right)\sum_{n\geq 0}\left\langle{n+k\atop k}\right\rangle_{d,F}x^n=N(d,k,x)$ *is a polynomial of degree* $(d-1)(k-1).$

*For* $d=2$ *we get more precisely*

$$\left(\sum_{j=0}^{2k+1}(-1)^{\binom{j+1}{2}}\binom{2k+1}{j}_F x^j\right)\sum_{n\geq 0}\left\langle{n+k\atop k}\right\rangle_{2,F}x^n=\sum_j\left\langle{k-1\atop j}\right\rangle_{2,F}x^j.$$

**Remark**

The polynomials $N(d,k,x)$ in general are not symmetric or unimodal. For example for $d=3$ and $k=5$ Mathematica gives

$N(3,5,x)=1+105x+9450x^2-7917x^3+166712x^4+7917x^5+9450x^6-105x^7+x^8.$

## 6. A common generalization

A generalization which contains all above cases is given by the Fibonacci polynomials

$$F_n(s,t)=\sum_{j=0}^{\left\lfloor\frac{n-1}{2}\right\rfloor}\binom{n-1-j}{j}t^j s^{n-2j}. \tag{42}$$

They satisfy $F_n(s,t)=sF_{n-1}(s,t)+tF_{n-2}(s,t)$ with initial values $F_0(s,t)=0$ and $F_1(s,t)=1.$

For $s=2$ and $t=-1$ we have $F_n(2,-1)=n$ which gives the original Hoggatt matrices, for $s=1+q$ and $t=-q$ we get $F_n(1+q,-q)=[n]_q$ which gives the $q-$analogs and for $s=t=1$ we get the Fibonacci analogs.

Binet's formulae give $F_n(s,t)=\dfrac{\alpha^n-\beta^n}{\alpha-\beta}$ with $\alpha=\dfrac{s+\sqrt{s^2+4t}}{2}$ and $\beta=\dfrac{s-\sqrt{s^2+4t}}{2}.$

We also have

$$\begin{pmatrix}0 & 1\\ t & s\end{pmatrix}^n=\begin{pmatrix}tF_{n-1}(s,t) & F_n(s,t)\\ tF_n(s,t) & F_{n+1}(s,t)\end{pmatrix},$$

From

$$\begin{pmatrix}tF_{n-1}(s,t) & F_n(s,t)\\ tF_n(s,t) & F_{n+1}(s,t)\end{pmatrix}=\begin{pmatrix}0 & 1\\ t & s\end{pmatrix}^{n-k}\begin{pmatrix}0 & 1\\ t & s\end{pmatrix}^k=\begin{pmatrix}tF_{n-k-1}(s,t) & F_{n-k}(s,t)\\ tF_{n-k}(s,t) & F_{n-k+1}(s,t)\end{pmatrix}\begin{pmatrix}tF_{k-1}(s,t) & F_k(s,t)\\ tF_k(s,t) & F_{k+1}(s,t)\end{pmatrix}$$

we get by comparing the top right elements

$$F_n(s,t)=F_{k+1}(s,t)F_{n-k}(s,t)+tF_{n-k-1}(s,t)F_k(s,t).$$

This is equivalent with

$$\binom{n}{k}_{F(s,t)}=F_{k+1}(s,t)\binom{n-1}{k}_{F(s,t)}+tF_{n-k-1}(s,t)\binom{n-1}{k-1}_{F(s,t)}, \tag{43}$$

which shows that the Fibonomials

$$\binom{n}{k}_{F(s,t)}=\prod_{j=0}^{k-1}\frac{F_{n-j}(s,t)}{F_{k-j}(s,t)} \tag{44}$$

are polynomials in $s,t$ with nonnegative integer coefficients.

A combinatorial proof of this fact has been given in [2] and an arithmetic one in [10].

The first terms of the Fibonomial triangle are

$$\begin{pmatrix} 1 & & & & & \\ 1 & 1 & & & & \\ 1 & s & 1 & & & \\ 1 & s^2+t & s^2+t & 1 & & \\ 1 & s\left(s^2+2t\right) & \left(s^2+t\right)\left(s^2+2t\right) & s\left(s^2+2t\right) & 1 & \\ 1 & s^4+3s^2t+t^2 & \left(s^2+2t\right)\left(s^4+3s^2t+t^2\right) & \left(s^2+2t\right)\left(s^4+3s^2t+t^2\right) & s^4+3s^2t+t^2 & 1 \end{pmatrix}$$

As above we get

$$h_k(x,s,t)=\sum_{j=0}^{k}(-1)^{\binom{j+1}{2}}t^{\binom{j}{2}}\binom{k}{j}_{F(s,t)}x^j=\prod_{j=0}^{\left\lfloor\frac{k}{2}\right\rfloor}u_{k-2j}\left((-t)^j x,s,t\right) \tag{45}$$

with $u_n(x,s,t)=1-L_{n-1}(s,t)x+(-t)^{n-1}x^2$ for $n>1$ and $u_0(x,s,t)=1$ and $u_1(x,s,t)=1-x,$

where the Lucas polynomials $L_n(s,t)$ satisfy $L_n(s,t)=sL_{n-1}(s,t)+tL_{n-2}(s,t)$ with initial values $L_0(s,t)=2$ and $L_1(s,t)=s.$

This implies as before that

$$\frac{1}{h_{k+1}(x,s,t)}=\frac{1}{\sum_{j=0}^{k+1}(-1)^{\binom{j+1}{2}}t^{\binom{j}{2}}\binom{k+1}{j}_{F(s,t)}x^j}=\sum_{n\geq 0}\binom{n+k}{k}_{F(s,t)}x^n. \tag{46}$$

The Hoggatt coefficients can be defined by

$$\left\langle\begin{matrix} n \\ k \end{matrix}\right\rangle_{d,F(s,t)}=\prod_{j=0}^{k-1}\frac{\binom{n-j+d-1}{d}_{F(s,t)}}{\binom{k-j+d-1}{d}_{F(s,t)}}. \tag{47}$$

An extension of a result which in [5] has been obtained for $d=2$ is

**Theorem 15**

$$\det\left((-t)^{\binom{i-j}{2}}\binom{n}{k-j+i}_{F(s,t)}\right)_{i,j=0}^{d-1}=\left\langle\begin{matrix} n \\ k \end{matrix}\right\rangle_{d,F(s,t)}. \tag{48}$$

This implies that all Hoggatt coefficients are polynomials in $s,t$ with integer coefficients.

**Proof**

By Binet's formula we have $F_n(s,t)=\alpha^{n-1}\dfrac{1-\left(\frac{\beta}{\alpha}\right)^n}{1-\left(\frac{\beta}{\alpha}\right)}$ with $\dfrac{\beta}{\alpha}=-\dfrac{t}{\alpha^2}$. If we set $\dfrac{\beta}{\alpha}=q,$ then

we get

$$F_n=\left(-\frac{q}{t}\right)^{\frac{1-n}{2}}\frac{1-q^n}{1-q}. \tag{49}$$

This implies

$$\binom{n}{k}_{F(s,t)}=\left(-\frac{q}{t}\right)^{\frac{k^2-nk}{2}}\begin{bmatrix} n \\ k \end{bmatrix}_q. \tag{50}$$

We can now formulate $\det\left((-t)^{\binom{i-j}{2}}\binom{n}{k-j+i}_{F(s,t)}\right)_{i,j=0}^{d-1}$ in terms of $q-$binomial coefficients as

$$\det\left((-t)^{\binom{i-j}{2}}\binom{n}{k-j+i}_{F(s,t)}\right)_{i,j=0}^{d-1}=\det\left((-t)^{\binom{i-j}{2}}\left(-\frac{q}{t}\right)^{\frac{(k-j+i)(k-j+i-n)}{2}}\begin{bmatrix} n \\ k-j+i \end{bmatrix}_q\right)_{i,j=0}^{d-1}$$

$$=(-1)^{\frac{k(k-n)d}{2}}t^{\frac{k(n-k)d}{2}}q^{\frac{k(k-n)d}{2}}\det\left(q^{\binom{i-j}{2}}\begin{bmatrix} n \\ k-j+i \end{bmatrix}_q\right)_{i,j=0}^{d-1}.$$

The last determinant has been computed above as

$$\left\langle{n \atop k}\right\rangle_{d,q}=\prod_{j=0}^{d-1}\frac{\begin{bmatrix} n+j \\ k \end{bmatrix}_q}{\begin{bmatrix} k+j \\ k \end{bmatrix}_q}=\prod_{j=0}^{d-1}\frac{\left(-\frac{t}{q}\right)^{\frac{k^2-(n+j)k}{2}}}{\left(-\frac{t}{q}\right)^{\frac{k^2-(k+j)k}{2}}}\prod_{j=0}^{d-1}\frac{\binom{n+j}{k}_{F(s,t)}}{\binom{n+j}{k}_{F(s,t)}}=\left(-\frac{t}{q}\right)^{\frac{-k(n-k)d}{2}}\left\langle{n \atop k}\right\rangle_{d,F(s,t)}.$$

This gives

$$\det\left((-t)^{\binom{i-j}{2}}\binom{n}{k-j+i}_{F(s,t)}\right)_{i,j=0}^{d-1}=(-1)^{\frac{k(k-n)d}{2}}t^{\frac{k(n-k)d}{2}}q^{\frac{k(k-n)d}{2}}\left(-\frac{t}{q}\right)^{\frac{-k(n-k)d}{2}}\left\langle{n \atop k}\right\rangle_{d,F(s,t)}=\left\langle{n \atop k}\right\rangle_{d,F(s,t)}.$$

The same proof as above gives

**Theorem 16**

$$\det\left(\binom{n+d+j-i-1}{n-1}_{F(s,t)}\right)_{i,j=0}^{k-1}=(-t)^{\binom{k}{2}d}\left\langle{n \atop k}\right\rangle_{d,F(s,t)}. \tag{51}$$

Let us also mention

**Theorem 17**

$$\frac{\det\left(\binom{n+i+j}{k+j}_{F(s,t)}\right)_{i,j=0}^{d-1}}{\det\left(\binom{n+i+j}{j}_{F(s,t)}\right)_{i,j=0}^{d-1}}=\left\langle{n \atop k}\right\rangle_{d,F(s,t)} \tag{52}$$

and

$$\frac{\det\left(\binom{n+i+k}{k+j}_{F(s,t)}\right)_{i,j=0}^{d-1}}{\det\left(\binom{n+i}{j}_{F(s,t)}\right)_{i,j=0}^{d-1}}=\left\langle{n+k \atop k}\right\rangle_{d,F(s,t)}. \tag{53}$$

The proof follows in the same way as in Theorem 11.

Let us consider two extreme special cases.

Taking limits for $t \to 0$ we get $F_n(s,0)=s^{n-1}$ for $n \geq 1$ and $F_0(s,0)=0.$

For the Fibonomials we get $\binom{n}{k}_{F(s,0)}=s^{k(n-k)}$ for $0 \leq k \leq n.$

This follows by induction from (43). From (47) we get that the entries of the Hoggatt matrices are

$$\left\langle{n \atop k}\right\rangle_{d,F(s,0)}=s^{dk(n-k)}.$$

Taking limits for $s \to 0$ gives more interesting results.

$F_{2n}(0,t)=0$ and $F_{2n+1}(0,t)=t^n$ by the definition of the Fibonacci polynomials.

The Lucas polynomials reduce to $L_{2n}(0,t)=2t^n$ and $L_{2n+1}(0,t)=0.$

Therefore we get

$$h_{2n}(x,0,t)=\left(t^{2n-1}x^2-1\right)^n \tag{54}$$

and

$$h_{2n+1}(x,0,t)=\left(1-t^n x\right)^{n+1}\left(1+t^n x\right)^n. \tag{55}$$

Comparing with (45) we get

$$\begin{aligned}&\binom{2n}{2j}_{F(0,t)}=\binom{n}{j}t^{2j(n-j)},\quad \binom{2n}{2j+1}_{F(0,t)}=0,\\ &\binom{2n+1}{2j}_{F(0,t)}=\binom{n}{j}t^{j(2n+1-2j)},\quad \binom{2n+1}{2j+1}_{F(0,t)}=\binom{n}{j}t^{(2n-1-2j)j+n}.\end{aligned} \tag{56}$$

For example

$$\left(\binom{n}{k}_{F(0,t)}\right)_{n,k=0}^{7}=\begin{pmatrix}1&&&&&&\\1&1&&&&&\\1&0&1&&&&\\1&t&t&1&&&\\1&0&2t^2&0&1&&\\1&t^2&2t^3&2t^3&t^2&1&\\1&0&3t^4&0&3t^4&0&1\end{pmatrix}. \tag{57}$$

Let us also mention the Hoggatt triangle $H_{2,F(0,t)}$.

Here we get by (56)

$$\left\langle {n \atop k} \right\rangle_{2,F(0,t)}=\prod_{j=0}^{k-1}\frac{\binom{n+1-j}{2}_{F(0,t)}}{\binom{k+1-j}{2}_{F(0,t)}}=\prod_{j=0}^{k-1}\frac{\left\lfloor\frac{n+1-j}{2}\right\rfloor t^{n-j}}{\left\lfloor\frac{k+1-j}{2}\right\rfloor t^{k-j}}=t^{k(n-k)}\binom{\left\lfloor\frac{n}{2}\right\rfloor}{\left\lfloor\frac{k}{2}\right\rfloor}\binom{\left\lfloor\frac{n+1}{2}\right\rfloor}{\left\lfloor\frac{k+1}{2}\right\rfloor}. \tag{58}$$

for $0\le k\le n$ and $=0$ else.

The first terms are

$$\left(\left\langle {n \atop k} \right\rangle_{2,F(0,t)}\right)_{n,k=0}^{7}=\begin{pmatrix}1&&&&&&\\1&1&&&&&\\1&t&1&&&&\\1&2t^2&2t^2&1&&&\\1&2t^3&4t^4&2t^3&1&&\\1&3t^4&6t^6&6t^6&3t^4&1&\\1&3t^5&9t^8&9t^9&9t^8&3t^5&1\end{pmatrix}. \tag{59}$$

For $t=1$ this is OEIS [8], A088855. For $t=1$ the row sums are $\binom{n+1}{\left\lfloor \frac{n+1}{2} \right\rfloor}$ and for $t=-1$ the sum of row $2n$ is the Catalan number $C_n$ and the sum of row $2n-1$ is the central binomial coefficient $\binom{2n}{n}$.

A companion to Conjecture 13 is

**Conjecture 18**

$$\left(\sum_{j=0}^{dk+1}(-1)^{\binom{j+1}{2}}t^{\binom{j}{2}}\binom{dk+1}{j}_{F(s,t)}x^j\right)\sum_{n\geq 0}\left\langle {n+k \atop k} \right\rangle_{d,F(s,t)} x^n = N(d,k,s,t,x) \quad (60)$$

*is a polynomial of degree* $(d-1)(k-1)$.

*For* $d=2$ *we get more precisely*

$$\sum_{j=0}^{2k+1}(-1)^{\binom{j+1}{2}}\binom{2k+1}{j}_F t^{\binom{j}{2}}x^j\sum_{n\geq 0}\left\langle {n+k \atop k} \right\rangle_{2,F(s,t)} x^n = \sum_j \left\langle {k-1 \atop j} \right\rangle_{2,F(s,t)} t^{j^2+j}x^j. \quad (61)$$

For $s\to 0$ identity (61) reduces to

$$\left(1-t^kx\right)^{k+1}\left(1+t^kx\right)^k\sum_{n=0}^{\infty}\binom{\left\lfloor \frac{n+k}{2}\right\rfloor}{\left\lfloor \frac{k}{2}\right\rfloor}\binom{\left\lfloor \frac{n+k+1}{2}\right\rfloor}{\left\lfloor \frac{k+1}{2}\right\rfloor}t^{kn}x^n = \sum_{j=0}^{k}\binom{\left\lfloor \frac{k-1}{2}\right\rfloor}{\left\lfloor \frac{j}{2}\right\rfloor}\binom{\left\lfloor \frac{k}{2}\right\rfloor}{\left\lfloor \frac{j}{2}\right\rfloor}t^{kj}x^j. \quad (62)$$

Let us prove this identity. It suffices to consider $t=1$.

Let $k=2\ell$. Then

$$\sum_{n=0}^{\infty}\binom{\left\lfloor \frac{2n+k}{2}\right\rfloor}{\left\lfloor \frac{k}{2}\right\rfloor}\binom{\left\lfloor \frac{2n+k+1}{2}\right\rfloor}{\left\lfloor \frac{k+1}{2}\right\rfloor}x^{2n} = \sum_{n=0}^{\infty}\binom{n+\ell}{\ell}^2 x^{2n},$$

$$\sum_{n=0}^{\infty}\binom{\left\lfloor \frac{2n+1+k}{2}\right\rfloor}{\left\lfloor \frac{k}{2}\right\rfloor}\binom{\left\lfloor \frac{2n+1+k+1}{2}\right\rfloor}{\left\lfloor \frac{k+1}{2}\right\rfloor}x^{2n+1} = \sum_{n=0}^{\infty}\binom{n+\ell}{\ell}\binom{n+1+\ell}{\ell}x^{2n+1}.$$

Thus the left-hand side is

$$(1-x)(1-x^2)^{2\ell}\left(\sum_{n=0}^{\infty}\binom{n+\ell}{\ell}^2 x^{2n}+\sum_{n=0}^{\infty}\binom{n+\ell}{\ell}\binom{n+1+\ell}{\ell}x^{2n+1}\right)$$

$$=(1-x)(1-x^2)^{2\ell}\left({}_2F_1\left[\begin{matrix}\ell+1,\ell+1\\1\end{matrix};x^2\right]+(\ell+1)x\ {}_2F_1\left[\begin{matrix}\ell+1,\ell+2\\2\end{matrix};x^2\right]\right).$$

The right-hand side reduces in an analogous way to ${}_2F_1\left[\begin{matrix}-\ell,-\ell+1\\1\end{matrix};x^2\right]+\ell x\ {}_2F_1\left[\begin{matrix}-\ell+1,-\ell+1\\2\end{matrix};x^2\right].$

Euler's transformation formula (cf. [1], (2.2.7))

${}_2F_1\left[\begin{matrix}a,b\\c\end{matrix};z\right]=(1-z)^{c-a-b}{}_2F_1\left[\begin{matrix}c-a,c-b\\c\end{matrix};z\right]$ gives

$${}_2F_1\left[\begin{matrix}\ell+1,\ell\\1\end{matrix};x^2\right](1-x^2)^{2\ell}={}_2F_1\left[\begin{matrix}-\ell,1-\ell\\1\end{matrix};x^2\right],$$

$${}_2F_1\left[\begin{matrix}\ell+1,\ell+1\\2\end{matrix};x^2\right](1-x^2)^{2\ell}={}_2F_1\left[\begin{matrix}1-\ell,1-\ell\\2\end{matrix};x^2\right].$$

By comparing coefficients we get

$${}_2F_1\left[\begin{matrix}\ell+1,\ell\\1\end{matrix};x^2\right]+\ell x\ {}_2F_1\left[\begin{matrix}\ell+1,\ell+1\\2\end{matrix};x^2\right]=(1-x)\ {}_2F_1\left[\begin{matrix}\ell+1,\ell+1\\1\end{matrix};x^2\right]+(\ell+1)x(1-x)\ {}_2F_1\left[\begin{matrix}\ell+1,\ell+2\\2\end{matrix};x^2\right]$$

which proves (62) for even $k$. In a similar way the formula is proved for odd $k$.